\documentclass[12pt]{amsart}
\usepackage{doublespace}

\usepackage{amsmath}
\usepackage{amsfonts}
\usepackage{latexsym}
\usepackage{amssymb}
\usepackage{amscd}

\newcommand{\Kt}{\widetilde{K}}
\newcommand{\Ft}{\widetilde{F}}
\newcommand{\Kc}{K^{cyc}}
\newcommand{\Fc}{F^{cyc}}
\newcommand{\Ki}{K_\infty}
\newcommand{\q}{\mathbb{Q}}
\newcommand{\G}{\text{Gal}}

\newcommand{\tor}{\operatorname{Tor}}
\newcommand{\g}{\mathcal{G}}

\newcommand{\oh}{\mathcal{O}}
\newcommand{\om}{\Omega}
\newcommand{\zp}{\mathbb{Z}_p}
\newcommand{\rt}{\rightarrow}
\newcommand{\lrt}{\longrightarrow}
\newcommand{\z}{\mathbb{Z}}
\newcommand{\la}{\Lambda}
\newcommand{\ex}{\operatorname{Ext}}
\newcommand{\hm}{\operatorname{Hom}}

\newtheorem{thm}{Theorem}
\newtheorem{lem}{Lemma}
\newtheorem{cor}{Corollary}
\newtheorem{con}{Conjecture}

\begin{document}

\title{Greenberg's Conjecture and Cyclotomic Towers}
\author{David C. Marshall}
\address{Department of Mathematics, University of Texas at Austin,
Austin, Texas 78712, USA}
\email{marshall@math.utexas.edu}
\subjclass[2000]{Primary 11R23; Secondary 11R18, 11R32}
\date{June 12, 2003}

\begin{abstract}
We describe Greenberg's pseudo-null conjecture, and prove a result describing conditions under which the pseudo-null conjecture for a number field $K$ implies the conjecture for finite extensions of $K$. We then apply the result to the cyclotomic $\zp$-tower above a cyclotomic field of prime roots of unity, verifying the conjecture for a large class of cyclotomic fields.  
\end{abstract}

\maketitle

\section{Greenberg's conjecture}
In the late 1950's Iwasawa introduced a powerful technique for studying class groups and unit groups of number fields. Motivated by the theory of curves over finite fields, Iwasawa's theory of $\zp$-extensions has since become a widely used tool in algebraic number theory, Galois theory, and arithmetic geometry.  We describe in this section a conjecture of Greenberg concerning the structure of a classical Iwasawa module, and we mention a Galois theoretic consequence concerning free pro-$p$-extensions of number fields.

Let $K$ be an algebraic number field and $p$ an odd prime. By a \textit{multiple $\zp$-extension} $\Ki/K$ we mean a Galois extension with Galois group $\Gamma\simeq \zp^d$ for some positive integer $d$. In what follows we will be particularly interested in two such extensions of $K$ for which we reserve the following notation:
\begin{itemize}
\item $\Kc/K$ denotes the \textit{cyclotomic} $\zp$-extension of $K$.
\item $\Kt/K$ denotes the compositum of all $\zp$-extensions of $K$.
\end{itemize}
Let $F$ be a finite extension of $K$ contained in $\Ki$, and denote by $A(F)$ the Sylow $p$-subgroup of the ideal class group of $F$. The Galois group of $F/K$ acts on $A(F)$ in the natural way, making $A(F)$ into a $\zp[\G(F/K)]$-module. As $F$ varies over all finite subextensions the $A(F)$ form an inverse system (under norm maps) and we denote by $A$ the inverse limit. The group $A$ then carries a natural structure as a module over the Iwasawa algebra 
\[\zp[[\Gamma]]:=\varprojlim_F\zp[\G(F/K)].\]

It is common to study $A$ by identifying the $A(F)$ with Galois groups as follows. By class field theory, the group $A(F)$ is isomorphic to the Galois group, $X_F$, of the maximal abelian unramified $p$-extension of $F$ (the \textit{$p$-Hilbert class field of $F$}). The isomorphism respects the Galois module structure, the action of $\G(F/K)$ on $X_F$ being inner automorphism.  The $X_F$ form an inverse system (the maps being given by restriction of automorphisms) and the limit $X$ is the Galois group of the maximal abelian unramified pro-$p$-extension of $\Ki$. So $X\simeq A$.

The Iwasawa algebra $\zp[[\Gamma]]$ is non-canonically isomorphic to the power series ring
\[\la:=\zp[[T_1, T_2, \dots, T_d]],\]
where topological generators $\gamma_i$ of $\Gamma$ are sent to $1+T_i$. So the $\zp[[\Gamma]]$-module structure of $A$ is studied via the $\la$-module structure of $X$ (noting that $T_ix=x^{\gamma_i-1}$).

For $\Ki/K$ any multiple $\zp$-extension Greenberg (\cite{Green:73}, Theorem~1) has shown $X$ to be a finitely generated torsion $\la$-module. In particular, the annihilator of $X$, $\operatorname{Ann}_{\la}(X)$, is non-trivial. Traditionally, annihilators of classical Iwasawa modules have been of much interest. The Main conjecture of Iwasawa theory gives the factors of the annihilator of $X$ for the cyclotomic $\zp$-extension of a number field $K$ as essentially the $p$-adic $L$-functions attached to $K$. There is also a two variable Main conjecture for certain $\zp^2$-extensions arising from the theory of elliptic curves.

Greenberg (\cite{Green:99}, Conjecture~3.4) has conjectured that for the cyclotomic $\zp$-extension $\Kc/K$ of a totally real field $K$, the module $X$ is finite. If a totally real field $K$ satisfies Leopoldt's conjecture the extensions $\Kc$ and $\Kt$ coincide (i.e. $K$ has only one $\zp$-extension). Furthermore, when $\la =\zp[[T]]$ it can be shown that a module being finite is equivalent to having an annihilator of height at least 2. With this in mind the above conjecture is a special case of the more general conjecture (\cite{Green:99}, Conjecture~3.5):

\begin{con}
Let $K$ be any number field and $\Kt$ the compositum of all $\zp$-extensions of $K$. Then $\operatorname{Ann}_{\la}(X)$ has height at least 2.
\end{con}

A $\la$-module whose annihilator has height at least 2 is said to be \textit{pseudo-null}, and we will refer to Conjecture~1 above as \textit{Greenberg's conjecture}, or just the \textit{pseudo-null conjecture}.

The point of this note is two-fold. First, we prove a ``going-up'' theorem for the pseudo-null conjecture. Namely, if $K$ is a number field, and $F$ 
is a finite extension of $K$ in $\Kt$, we give conditions under which Greenberg's conjecture for $K$ implies Greenberg's conjecture for $F$ (Theorem~6).
The result is an exercise in utilizing several equivalent
formulations of the conjecture. Versions of these formulations have appeared in Lannuzel and Nguyen-Quang-Do (\cite{Lan:00}, Theorem~4.4) 
as well as work of McCallum~\cite{McCal:00} and this author~\cite{Ma:00}.
Secondly, as an application of the result, we consider the example $K=\q(\zeta_p)$ and $F=\q(\zeta_{p^n})$. We verify the conjecture for a certain class of such $K$'s, implying the conjecture for each field in the corresponding $\zp$-tower. 

The key argument in both results is reduced to a capitulation problem, namely the need for a set of ideals, or ideal classes, to become principal when extended to an appropriate field. For the ``going-up'' result, the resolution of this problem is provided by an equivalent form of the conjecture, stating that all ideal classes capitulate in $\Kt$. In verifying the conjecture for $\q(\zeta_p)$ capitulation is obtained by
more direct means. We state our second result here.

Let $K=\q(\zeta_p)$, $E=\oh_K^{\times}$ and $U=\oh_{K_\pi}^{\times}$, where $\pi$ is the unique prime of $K$ above $p$.  
Denote by $\overline{E}$ the closure of $E$ in $U$. We denote by $\lambda_p$ the Iwasawa lambda invariant of the cyclotomic $\zp$-extension of $\q(\zeta_p)$. Let $v_p$ denote the $p$-adic valuation. In Section~4 we prove
\begin{thm}
Suppose $K=\q(\zeta_p)$ satisfies the following conditions:
\begin{enumerate}
\item  Vandiver's conjecture

\item $\lambda_p=1$.

\item $v_p(|(U/\overline{E})[p^\infty]|)\leq v_p(|A(K)|)$.
\end{enumerate}
Then  for all $n\geq 1$ the pseudo-null conjecture holds for $\q(\zeta_{p^n})$.
\end{thm}

We mention here one Galois theoretic consequence of the pseudo-null conjecture for cyclotomic fields. The existence of free pro-$p$-extensions (Galois extensions with Galois group a free pro-$p$-group) has been the subject of much study. See for example the list of known results in~\cite{Yama:94}. Let $K=\q(\zeta_{p^n})$ for some $n>0$, and let $\om_K$ denote the maximal pro-$p$ extension of $K$ which is unramified at all primes not dividing $p$. Let $\g_K$ denote the Galois group.

Since free pro-$p$-extensions are unramified outside $p$, such extensions of $K$ are contained in $\om_K$. We will see that $\g_K$ is
a free pro-$p$ group exactly when $p$ is a regular prime (since the
number of relations defining $\g_K$ is equal to the $p$-rank of the class group of $K$). When $p$ is an irregular prime the group $\g_K$ is not free, but we
may look for free pro-$p$ quotients. 
Let $r_2$ denote the number of complex places of $K$.
Then Leopoldt's conjecture predicts $r_2+1$ independent $\zp$-extensions
of $K$, and so the maximal rank of a free pro-$p$-extension of $K$
is bounded above by $r_2+1$. The following is proved in~\cite{Lan:00}, as well as~\cite{McCal:00}:
\begin{thm}
Suppose that $K=\q(\zeta_{p^n})$ satisfies Greenberg's conjecture. Then $\g_K$ has a free pro-$p$-quotient of rank $r_2+1$ if and only if $p$ is regular.
\end{thm}

We give here a brief outline of the paper. In Section 2, we introduce several auxiliary $\la$-modules and Galois  
groups needed for the later study. Theorem~3 and Lemma~1 are the key results of this section, implying a sufficient condition for a standard Iwasawa module to be torsion free (Corollary~1). In Section 3 we
recall and provide several equivalent formulations of Greenberg's pseudo-null conjecture, and we state and prove one of
our main results (the ``going-up'' theorem). Finally, in Section 4 we turn to the example furnished by cyclotomic fields, proving Theorem~1 above.

\textit{Acknowledgements} This work is an outgrowth of the authors Ph.D. thesis and he would like to thank his advisor Bill McCallum, as well as Ralph Greenberg and Manfred Kolster for useful conversations and comments. This work was partially supported by NSF VIGRE grant 9977116. 

\section{Auxiliary modules}
For a number field $K$ and a prime number $p$, we call a field extension of $K$ \textit{$p$-ramified} if it is unramified at all primes of $K$ not dividing $p$. We fix the following notation:

\noindent The fields:
\vspace{.1in}

\begin{tabular}{lll}
$\om_K$ & & the maximal pro-$p$, $p$-ramified extension of $K$ \\
$\Kt$ &   & the compositum of all $\zp$-extensions of $K$ \\
$L_\infty$ &  & the maximal abelian unramified pro-$p$-extension of $\Kt$ \\
$M_\infty$ &  & the maximal abelian $p$-ramified pro-$p$-extension of $\Kt$ \\
$N_\infty$ &  & the extension of $\Kt$ generated by $p$-power roots of $p$-units of $\Kt$ 
\end{tabular}
\vspace{.1in}

\noindent The Galois groups:
\vspace{.1in}

\begin{tabular}{lll}
$\g_K$ &  & the Galois group of $\om_K/K$ \\
$\Gamma$ &  & the Galois group of $\Kt/K$ \\
$X$ &  & the Galois group of $L_\infty/\Kt$ \\
$Y$ &  & the Galois group of $M_\infty/\Kt$ \\
$Y'$ &  & the Galois group of $N_\infty/\Kt$
\end{tabular}
\vspace{.1in}

The Galois groups $Y$ and $Y'$ carry an action of $\Gamma$ via conjugation, just as $X$, making them into $\la$-modules. We shall see that for certain base fields $K$, the pseudo-null conjecture may be formulated in terms of the $\la$-module structure of $Y$ (in particular, that $Y$ is $\la$-torsion free).  
The module $Y$ is known to be finitely generated, and, for $K/\q$ abelian, have $\la$-rank equal to $r_2$, where $r_2$ denotes the number of complex places of $K$ (\cite{Green:78}). For a $\la$-module $M$
 we write $\operatorname{Tor}_\la(M)$ for the $\la$-torsion submodule.
The following result  is due to
McCallum.
\begin{thm}[\cite{McCal:00}, Theorem 3]
Suppose there is only one prime of $K$ above $p$,
and $\Kt$ contains all $p$-power roots of unity. Then 
$\operatorname{Tor}_\la(Y')=0$.
\end{thm}

\textit{Remark 1}: The proof of this result involves a detailed analysis of the filtration 
\[E_F^u \subset E_F^n \subset E_F^{\text{loc}} \subset E_F,\]
where $E_F$ denotes the units $\oh_F[1/p]^\times$ of a finite extension $F$ of $K$ in $\Kt$, and the superscripts denote certain classes of universal norms (see Section~4 of~\cite{McCal:00} for the precise definitions).  The torsion submodule of $Y'$ is contained in the kernel of a surjective map of Galois groups. The Pontryagin dual of this kernel is $\varinjlim_F(E_F/E_F^u) \otimes \q_p/\z_p$, and is shown to be zero by considering each graded factor from the filtration.

\textit{Remark 2}: In particular, the result tells us $\tor_{\la}(Y)$ fixes the field $N_\infty$. This observation, combined with Lemma~1 below, gives our approach to verifying the pseudo-null conjecture.

The group $\g_K$ has a minimal free presentation
\[1\lrt R\lrt F_g\lrt \g_K \lrt 1,\]
where $F_g$ is the free pro-$p$-group on $g$ generators and $R$ is the normal closure of a finitely generated subgroup (the group of relations for $\g_K$). Denote by $s$ the minimal number of (topological) generators of $R$. The numbers $g$ and $s$ are equal to the $\mathbb{F}_p$-dimensions of $H^i(\g_K, \z/p\z)$, $i=1, 2$ respectively (see Chapter~4 of~\cite{Serre:97}).

Let $\g_K^{ab}$ denote the maximal abelian quotient of $\g_K$, and $M_K$ the maximal abelian $p$-ramified pro-$p$-extension of $K$ (so $\g_K^{ab}=\G(M_K/K)$).
The field $M_K$ is an abelian, $p$-ramified extension of $\Kt$ (the Galois group of $M_K/\Kt$ is just the torsion subgroup of $\g_K^{ab}$), and so is contained in the field $M_\infty$. Hence we have a natural map from $Y$ to $\g_K^{ab}$ given by restriction of automorphisms. We refer the reader to \cite{McCal:00} for a proof of the following.

\begin{lem}[\cite{McCal:00}, Lemma 24]
Suppose $K$ satisfies Leopoldt's conjecture. If $\g_K$ is a one-relator group (i.e. $s=1$), then the map
\[\operatorname{Tor}_\la(Y)\lrt \g_K^{ab}\]
is the zero map if and only if $\operatorname{Tor}_\la(Y)=0$.
\end{lem}

The following is an immediate consequence of Theorem~3 and Lemma~1:
\begin{cor}
If $K$ 
is a number field satisfying the hypotheses of Theorem~3 and Lemma~1,
then 
\begin{equation}
M_K\subset N_\infty \,\, \text{implies}\,\, \operatorname{Tor}_\la(Y)=0.
\end{equation}
\end{cor} 

\section{Equivalent formulations}
We have introduced the natural Iwasawa modules $X$ and $Y$ in the last section. The Galois action on each of the $X_F$ is also compatible with regard to extensions of ideal classes, so we may form the $\la$-module $\varinjlim_FX_F$ as well. Recall the groups $\ex^i_\la(\cdot, \la)$ are the right derived functors of $\hm_\la(\cdot, \la)$.

\begin{thm}
Let $p$ be an odd prime and let $K$ be a number field with a unique prime above $p$. Then $\ex^1_\la(X, \la)$ is the Pontryagin dual of $\varinjlim_FX_F$, where the $F$ vary over the finite extensions of $K$ in $\Kt$.
\end{thm}

\textbf{Proof}: Let $\mathfrak{m}$ denote the unique maximal ideal of $\la=\zp[[T_1, \dots, T_r]]$, and define
\[\omega_n(T_i)=(1+T_i)^{p^n}-1.\]
The result is obtained by establishing the isomorphism
\begin{equation}
H_{\mathfrak{m}}^r(X)\simeq \varinjlim_FX_F,
\end{equation}
where $H_{\mathfrak{m}}^i(X)$ denotes Grothendieck's local cohomology relative to the $\mathfrak{m}$-primary sequences
\[\textbf{x}_n=(p^n, \omega_n(T_1), \dots, \omega_n(T_r)).\]
The desired result is then a consequence of (a version of) Grothendieck's local duality; namely
\[\ex_\la^{N-i}(X, \la)\simeq \hm_{\zp}(H_{\mathfrak{m}}^i(X), \q/\z),\]
where $N$ denotes the length of the $\mathfrak{m}$-primary sequence.
A good reference for this material is Chapter~3 of~\cite{Bruns:93}.

The details establishing (2) can be found in Theorem~8 of~\cite{McCal:00}, where McCallum proves a similar result for the Galois group $X'$ of the maximal abelian unramified pro-$p$-extension of $\Kt$ in which all primes dividing $p$ are completely decomposed. The proof translates easily to this case, simply replacing the decomposition group with inertia. $\Box$

Let $\mu_n$ denote the group of $n$-th roots of unity. As above, we let $X_F'$ denote the Galois group of the maximal abelian unramified extension of $F$ in which all primes dividing $p$ are completely decomposed. We write $X'$ for $X_{\Kt}'$. 

\begin{thm}
Let $p>5$ be a prime and suppose $\mu_p$ is in $K$. If $K$ has a unique prime ideal $\wp$ dividing $p$, 
then the following are equivalent:

(a) $X$ is pseudo-null

(b) $X'$ is pseudo-null

(c) $\operatorname{Tor}_\la(Y)=0$

(d) $\varinjlim_FX_F'=0$

(e) $\varinjlim_FX_F=0$,

\noindent where the fields $F$ vary over all finite extensions of $K$ in $\Kt$.
\end{thm}

\textbf{Proof}:
$(a)\Leftrightarrow (b)$. Recall $\Gamma =\G(\Kt/K)$. We let $\Gamma_\wp$ denote the decomposition group of $\wp$ in $\Gamma$, and let $\la_\wp=\zp[[\Gamma/\Gamma_\wp]]$. There is a natural surjection $X\rightarrow X'$ whose kernel is generated as a $\zp$-module by the Frobenius automorphisms corresponding to the primes above $p$, and therefore is finitely generated as a module over $\la_\wp$. As a $\la$-module, the annihilator of $\la_\wp$ has height equal to the $\zp$-rank of $\Gamma_\wp$ (this is just the augmentation ideal in $\zp[[\Gamma_p]]$). Since there is only one prime of $K$ above $p$, its decomposition group has finite index in $\Gamma$, and therefore our assumption on $p$ makes $\la_\wp$ pseudo-null. Therefore the kernel of the surjection $X\rightarrow X'$ is pseudo-null, and $X$ and $X'$ are pseudo-isomorphic. 

$(a)\Leftrightarrow (c)$.  This follows from a duality due to Jannsen (\cite{Jan:89}, Theorem~5.4) relating the $\la$-modules $X'$ and $Y$, together with a
structure theorem for $Y$ due to Nguyen-Quang-Do (Corollary~14 of~\cite{McCal:00} or Theorem~4.4 of~\cite{Lan:00}). 

$(c)\Leftrightarrow (d)$. In proving the results cited in the previous case, one shows, in particular, that
\[\tor_\la(Y)\simeq \ex^1_\la(X', \la)\]
(\cite{McCal:00}, Theorem~9).
But $\ex^1_\la(X', \la)$ is known to be the Pontryagin dual of $\varinjlim_FX_F'$ (\cite{McCal:00}, Theorem~8). The result then follows.

$(c)\Leftrightarrow (e)$. Grothendieck's local duality can be used to show that a torsion $\la$-module is pseudo-null if and only if $\ex_\la^1$ vanishes 
(\cite{McCal:00}, Lemma~6). This implies, in particular, that $\ex^1_\la(X, \la)$ and $\ex^1_\la(X', \la)$ are isomorphic, yielding
\[\tor_\la(Y)\simeq \ex^1_\la(X, \la)\]
as well. Theorem~4 then finishes the proof. $\Box$

\textit{Remark}: Various forms of these equivalences have certainly appeared elsewhere. In~\cite{Lan:00}, Lannuzel and Nguyen-Quang-Do prove the equivalence of (a), (c), and (e) under slightly different hypotheses. Namely, no restriction is made on the number of primes of $K$ dividing $p$, but rather it is assumed that all finite extensions of $K$ in $\Kt$ satisfy Leopoldt's conjecture. 
Formulation (c) has been used by McCallum~\cite{McCal:00} and this author~\cite{Ma:00} to verify Greenberg's conjecture for certain classes of cyclotomic fields. 

The following theorem provides sufficient conditions for when the pseudo-null conjecture for a number field $K$ implies the conjecture for a finite extension of $K$ in $\Kt$. We apply this to the cyclotomic tower in Section~4.

\begin{thm}
Let $p\geq 5$ be a prime and suppose $\mu_p$ is contained in $K$. Suppose $K$ has a unique prime 
$\wp$ dividing $p$. Then, if $F\subset \Kt$
is a finite extension of $K$ satisfying
\begin{enumerate}
\item $\wp$ is non-split in $F/K$

\item $\dim_{\mathbb{F}_p}H^2(\g_F, \z/p\z)\leq 1$

\item Leopoldt's conjecture,
\end{enumerate} 
then Greenberg's conjecture for $K$ implies Greenberg's conjecture
for $F$. 
\end{thm} 

\textbf{Proof}:
Let $K$ and $F$ be number fields satisfying the above hypotheses, and assume 
the pseudo-null conjecture holds for $K$. We apply the notation introduced 
in Section 2 to the field $F$ (so we have $\om_F$, $\g_F$, $M_F$, etc.) If the $\mathbb{F}_p$-dimension of $H^2(\g_F, \z/p\z)$ is 0, then $\g_F$ is a free pro-$p$-group. A structure theorem for $Y$ due to Nguyen Quang Do (\cite{NQD:84},
Proposition 1.7)
then implies $\tor_\la(Y)=0$. Hence by formulation (c) of Theorem~5 Greenberg's conjecture holds for $F$.

If the $\mathbb{F}_p$-dimension of $H^2(\g_F, \z/p\z)$ is 1, then
such an $F$ satisfies the hypotheses of Theorem~3 and Lemma~1, and so Corollary~1 applies.
Namely, Greenberg's pseudo-null conjecture will hold for $F$ 
provided $M_F\subset N_\infty$, and hence 
it will suffice to show the extension $M_F/\Ft$ is generated by 
$p$-power roots of $p$-units of $\Ft$.

We consider the field $\Fc = F\Kc$, the cyclotomic $\zp$-extension of $F$. By assumption, this field contains all $p$-power roots of unity. Recall the group $\g_F^{ab}=\G(M_F/F)$. The subgroup $\G(M_F/\Fc)$ has the same torsion subgroup (which is just $\G(M_F/\Ft)$) and $\zp$-rank 1 less. In particular, we have a non-canonical isomorphism
\[\G(M_F/\Fc)\simeq \G(\Ft/\Fc) \times \G(M_F/\Ft).\]
We let $L$ denote the fixed field of the first factor (so $M_F=\Ft L$.)

The Galois group $\G(L/\Fc)$ is isomorphic to the torsion subgroup of $\g_F^{ab}$, and hence is a finite $p$-group. Since $\Fc$ contains all $p$-power roots of unity, the extension $L/\Fc$ is just a Kummer extension, generated by $p$-power roots of elements of $\Fc$,
\[L=\Fc(x_1^{1/p^{m_1}}, x_2^{1/p^{m_2}}, \dots, 
x_n^{1/p^{m_n}}).\] 
Further, the ideals $(x_i)$ are $p^{m_i}$-th powers of ideals of $\Fc$, say $(x_i)=\mathfrak{J}_i^{p^{m_i}}$.

The extension $M_F/\Ft$ is also generated by the $x_i^{1/p^{m_i}}$, and 
the ideals $(x_i)$ are the $p^{m_i}$-th powers of the ideals 
$\mathfrak{J}_i$ extended to $\Ft$. But here is the key: the ideal 
classes $[\mathfrak{J}_i]$ become \textit{principal classes} when 
extended to $\Ft$. This follows from the fact that $\Fc \subset \Kt$ and, 
having assumed the pseudo-null conjecture holds for $K$ (using formulation 
(e) of Theorem~5), the fact that all ideal classes become principal in $\Kt$.

For a generator $x_i^{1/p^{m_i}}$ of $M_F/\Ft$  we now know the ideal 
$(x_i)$ is the $p^{m_i}$-th power of a principal ideal, say
\[(x_i)=(y_i)^{p^{m_i}}.\]
The elements $x_i$ and $y_i^{p^{m_i}}$ must differ by a unit, 
say $x_i=uy_i^{p^{m_i}}$. But clearly, an extension generated by a 
$p^{m_i}$-th root of $x_i$ is also generated by a $p^{m_i}$-th root  of 
$x_i/(y_i^{p^{m_i}})=u$, and so the extension $M_F/\Ft$ is generated by 
$p$-power roots of units on $\Ft$. This implies $M_F\subset N_\infty$ 
which, by Corollary~1 and Theorem~5, 
implies Greenberg's conjecture for $F$. $\Box$

\section{Cyclotomic Fields}
We fix $p$ a prime number and consider more closely the case of the cyclotomic fields $K=\q(\zeta_{p^n})$. Recall the group $\g_K$ has a minimal presentation as a pro-$p$-group with $g$ generators and $s$ relations, where $g$ and $s$ are equal to the $\mathbb{F}_p$-dimensions of $H^1(\g_K, \z/p\z)$ and $H^2(\g_K, \z/p\z)$ respectively. 
\begin{lem}
Let $p$ be a prime and let $K=\q(\zeta_{p^n})$ for some natural number $n$. Let $\alpha$ denote the $\z/p\z$-rank of the $p$-class group of $K$. Then 
\begin{align*}
g & = \frac{p^n + p^{n-1} +2}{2} +\alpha \\
s & = \alpha .
\end{align*}
\end{lem}

\textbf{Proof}: These computations are not new, and we give here just a sketch. 
Let $\Omega_K'$ be the maximal $p$-ramified extension of $K$ with Galois group $\g_K'$. Since $K$
contains the group $\mu_p$, and $\g_K$ is the maximal pro-$p$ quotient
of $\g_K'$, we have 
\[H^i(\g_K, \z/p\z)\simeq H^i(\g_K', \mu_p).\]
The $\z/p\z$-dimensions of the latter groups can be obtained by considering the sequence
\[1\lrt \mu_p \lrt \oh_{\om_K'}[1/p]^\times \stackrel{p}{\lrt}  
\oh_{\om_K'}[1/p]^\times \lrt  1.\] 
The $p$-power map on $\oh_{\om_K'}[1/p]^\times$ is surjective by the maximality 
of $\om_K'$ over $K$ (since $p$-th roots of $p$-units generate 
$p$-ramified extensions). 
Taking cohomology of the sequence with respect
to the Galois group $\g_K'$ yields a long exact sequence which may be broken
into the following pair of short exact sequences.
\begin{equation*}
0\rt \frac{\oh_K[1/p]^\times}{(\oh_K[1/p]^\times)^p} \rt H^1(\g_K', \mu_p) \rt
C(K)[p]\rt 0 
\end{equation*}

\begin{equation*}
0\rt \frac{C(K)}{pC(K)}
\rt H^2(\g_K', \mu_p)\rt H^2(\g_K', \oh_{\om_K'}[1/p]^\times)[p]\rt 0,
\end{equation*}
where $C(K)$ denotes the ideal class group of $K$.
The group $H^2(\g_K', \oh_{\om_K'}[1/p]^\times)$ injects into the Brauer group $B(K)$, and can be shown to be 0 by considering its behavior in the exact sequence
\begin{equation*}
0\rt B(K)\rt \oplus_vB(K_v)\stackrel{\sum inv}\lrt \q/\z \rt 0.
\end{equation*}
A simple dimension count then gives
\begin{align*}
g & = r_2+1+\alpha \\
s & = \alpha 
\end{align*}
where $r_2=(p^n +p^{n-1})/2$, as desired. $\Box$

If $p$ is a regular prime, $\alpha=0$ for $\q(\zeta_{p^n})$, $n\geq 0$. 
Hence $s=0$, implying $\tor_\la(Y)=0$, establishing Greenberg's conjecture for each field in the cyclotomic tower. 

The following corollary is an immediate consequence of Theorem~6 and Lemma~2.
\begin{cor}
Let $p$ be an irregular prime. Let $n>0$ be such that $\q(\zeta_{p^n})$ has a 
cyclic $p$-class group. Then Greenberg's conjecture for $\q(\zeta_p)$ implies Greenberg's conjecture for $\q(\zeta_{p^n})$.
\end{cor}

\textbf{Proof}: In the notation of Theorem~6, with $K$ as above, let $F=\q(\zeta_{p^n})$ for some positive integer $n$ satisfying the hypothesis. The field $K$ has a unique prime $\pi$ above $p$, and $\pi$ is totally ramified in $F/K$, and hence non-split. The dimension of $H^2(\g_F, \z/p\z)$ is less than or equal to 1 by our assumption of cyclic $p$-class groups.
Since $F/\q$ is abelian, implying Leopoldt's conjecture for $F$, the hypotheses of Theorem~6 are satisfied, as desired. $\Box$

Finally, we prove Theorem~1 by providing a class of cyclotomic fields $\q(\zeta_p)$, satisfying the hypotheses of Corollary~2, for which the pseudo-null conjecture is true. A similar class was first given by McCallum (\cite{McCal:00}, Theorem~1). He considered such fields with $p$-class group isomorphic to $\z/p\z$. We provide here a slight generalization of that class, allowing for cyclic $p$-class groups of arbitrary $p$-power order, as well as apply Corollary~2 to extend the conjecture to all fields in the cyclotomic $\zp$-tower. We restate Theorem~1 here.

\begin{thm}
Suppose $K=\q(\zeta_p)$ satisfies the following conditions:
\begin{enumerate}
\item Vandiver's conjecture

\item $\lambda_p=1$.

\item $v_p(|(U/\overline{E})[p^\infty]|)\leq v_p(|A(K)|)$.
\end{enumerate}
Then  for all $n\geq 1$ the pseudo-null conjecture holds for $\q(\zeta_{p^n})$.
\end{thm}

\textit{Remark 1}: Condition (2) is heuristically true for approximately
75\% of all irregular primes and experimentally true for 75\% of the 
irregular primes up to 12 million, according to~\cite{Buhler:93} 
(for these primes, $\lambda_p$ is just the index of irregularity of $p$). 

\textit{Remark 2}:  Letting $K_n=\q(\zeta_{p^{n+1}})$ and
$A_n=A(K_n)$, the hypotheses of Vandiver's conjecture and
$\lambda_p=1$ imply
\[A_n\simeq X/((1+T)^{p^n})-1)X,\]
where $X=\zp[[T]]/(T+p^a)$ (see Theorem 10.16 and Proposition 13.22
of~\cite{Wash:97}). In particular this yields isomorphisms
\[A_n\simeq \z/p^{a+n}\z\] for all $n\geq 0$, and so (3) is a 
condition on cyclic groups of $p$-power order.

\textit{Remark 3}: Since $A(K)$ is cyclic, there is only one Bernoulli 
number $B_i$, $2\leq i\leq p-3$, divisible by $p$. 
If $B_{p-j}$ denotes this 
term (so $\varepsilon_j A(K)$ is the non-trivial term of the 
idempotent decomposition of $A(K)$),then 
$L_p(s, \omega^{1-j})$ is the only non-trivial 
$p$-adic $L$-function attached to $K$.   
It follows from Theorem 8.25 of~\cite{Wash:97} that
\[(U/\overline{E})[p^{\infty}]\simeq \z/p^m\z ,\]
where $m=v_p(L_p(1, \omega^{1-j}))$.
This valuation may be computed in terms of the characteristic power 
series $f(T)$ of $\varprojlim_nA(K_n)$.
Under the assumption $\lambda_p=1$ this power series has the form 
$f(T)=(T+cp^a)u$, where $u$ is a unit, $p^a$ is the order of the 
cyclic group $A(K)$, and
\[f((1+p)^s-1)=L_p(s, \omega^{1-j}).\]
So the valuation of $L_p$ at $s=1$ equals the valuation of $f(p)=(p+cp^a)u$.

If $a>1$, $v_p(f(p))=1$, and condition (3) is satisfied. If, on the other 
hand, $a=1$, $v_p(f(p))$ depends on the value of $c\pmod{p}$. 
The valuation will again be 1 provided $c \not\equiv -1 \pmod{p}$. This 
congruence has been checked for $p<4000$ in~\cite{Iwa:65},
although tables are only given for $p<400$ and $3600<p<4000$. For these 
values the congruence condition is satisfied.

Suppose $K=\q(\zeta_p)$ satisfies (1)-(3) above. Since $A(K)$ is cyclic, say of order $p^a$, the group $\g$ is a one-relator group and Lemma~1 applies.
We will utilize this lemma to show $\tor_\la(Y)=0$. 
In light of Corollary~1, it suffices to show $M_K\subset N_\infty$, and so 
we consider the structure of $\g_K^{ab}$ in more detail.

\begin{lem}
Suppose $K$ satisfies hypothesis (1) and (2) of Theorem 7. Then the torsion
subgroup of $\g_K^{ab}$ is cyclic.
\end{lem}

\textbf{Proof}:
Let $J_K$ denote the idele group of $K$, with $K^\times$ embedded diagonally. 
Let $U$ be the subgroup of ideles which are units at $\pi$ (the prime of $K$ 
above $p$) and 1 elsewhere, 
and let $U'$ be the subgroup of ideles which are 1 at $\pi$ and units 
elsewhere. Class field theory gives an isomorphism
\[\g_K^{ab}\simeq \text{pro-$p$-completion of}\,\, 
J_K/(\overline{K^\times U'}),\]
where the overline denotes the closure.

If we let $\overline{E}$ denote the closure of the embedding of the units of $K$ in $U$, then in fact we have an exact sequence
\[0\lrt U_1/\overline{E}_1 \lrt \g_K^{ab} \lrt A(K) \lrt 0,\]
where the subscript 1 indicates we are taking units congruent 
to 1 modulo $\pi$. Since $U_1$ has $\zp$-rank
$[K:\q]=p-1$ and $\overline{E_1}$ has $\zp$-rank $(p-3)/2$ 
(by Leopoldt's conjecture, which holds for $K$), the $\zp$-rank of 
$\g_K^{ab}$ is $(p+1)/2$ ($p\neq 2$ by the assumption $\lambda_p =1$).

We claim the torsion in $\g_K^{ab}$ comes from $U_1/\overline{E}_1$,
and show this by considering an idele $(a_v)$ whose image in $\g_K^{ab}$
is a torsion element, say of order $p^m$. So
\[(a_v)^{p^m} \in \overline{K^\times U'},\]
say $(a_v)^{p^m}=\alpha (u_v)$ (where we abuse notation writing $\alpha$
for both the element of $K^\times$ as well as its diagonal image in $J_K$). 
This implies $\alpha$ is a $p^m$-th power in $K_\pi$, the 
$\pi$-adic completion of $K$. Let $\mathfrak{a}$ then be the ideal of $K$ 
such that $\mathfrak{a}^{p^m}=(\alpha)$. We want to show the class of
 $\mathfrak{a}$ is principal.

Let $K_{m-1}=\q(\zeta_{p^m})$, so 
$K_{m-1}(\alpha^{1/p^m})$ is an unramified extension. Since the 
class of $\mathfrak{a}$
lies in $A(K)^-$ (by Vandiver's conjecture), the Kummer pairing implies the Galois group of $K_{m-1}(\alpha^{1/p^m})/K_{m-1}$ is trivial. 
Hence $\alpha$ must be a $p^m$-th power in $K_{m-1}$ as well, which
means the ideal class of $\mathfrak{a}$ is principal when 
extended to $K_{m-1}$ (represented by a principal ideal generated by a 
$p^m$-th root of $\alpha$). 
But the map from $A(K)$ to $A(K_{m-1})$ is injective (\cite{Wash:97}, 
Proposition~13.26), and so $\mathfrak{a}$ must have represented a principal 
class in $A(K)$ as well. Hence the torsion in $\g_K^{ab}$ maps to 0 in $A(K)$.

We now just need to determine the torsion subgroup of 
$U_1/\overline{E}_1$. We may consider each factor of the idempotent decomposition separately. Since $\varepsilon_iE_1=0$ for $i=0$ and
for $i$ odd, and each  $\varepsilon_iU_1\simeq \zp$, we obtain
\[U_1/\overline{E}_1\simeq (\zp)^{(p+1)/2}\oplus \bigoplus_{\text{$i$ even}}
\varepsilon_iU_1/\varepsilon_i\overline{E_1}.\]
For even $i$ the terms $\varepsilon_iU_1/\varepsilon_i\overline{E_1}$ are
equal to $\varepsilon_iU_1^+/\varepsilon_i\overline{E_1}^+$, where the superscript $+$ indicates we are looking at units in the local subfield fixed by the automorphism of order 2. Vandiver's conjecture implies the cyclotomic
units $C_1^+$ have index prime to $p$ in $E_1^+$ (\cite{Wash:97}, Theorem~8.2), and so it suffices
to consider the quotients 
$\varepsilon_iU_1^+/\varepsilon_i\overline{C_1}^+$. But Theorem~8.25 of~\cite{Wash:97} states
\[[\varepsilon_iU_1^+:\varepsilon_i\overline{C_1}^+]=
p^{v_p(L_p(1, \omega^i))}.\]
Since $A(K)$ is cyclic there is only one non-trivial $L_p(s, \omega^i)$, and hence only one cyclic factor, say of order $p^m$, in the torsion subgroup of $U_1/\overline{E}_1$.  $\Box$

\textbf{Proof of Theorem 7}:
The field $\Kt$ is in fact the fixed field of the torsion subgroup of $\g_K^{ab}$, and so the extension $M_K/\Kt$ is a Kummer extension with $\G(M_K/\Kt)\simeq \z/p^m\z$. With $A(K)\simeq \z/p^a\z$, condition (3) of the Theorem just states $m\leq a$.

To show that $M_K$ is contained in $N_\infty$, we need to show
that $M_K/\Kt$ is generated by a $p$-th power root of a unit 
of $\Kt$. The argument, as in the proof of Theorem~6, is reduced to
a capitulation problem.

Consider the extension $M_K/K_{m-1}$. There is a non-canonical
isomorphism
\[\G(M_K/K_{m-1})\simeq \G(\Kt/K_{m-1})\times 
\G(M_K/\Kt).\]
We let $L$ denote the fixed field of the first factor.  
The extension $L/K_{m-1}$ is a Kummer extension, and we may write
\[L=K_{m-1}(x^{1/p^m})\] for some $x$ in $K_{m-1}$ where 
the ideal $(x)$ is of the form $(x)=\mathfrak{J}^{p^m}P$, where $P$ is 
the principal ideal of $K_{m-1}$ lying above $p$.

Since, in particular, $\mathfrak{J}$ represents a class of order dividing
$p^m$ in $A(K_{m-1})$, condition (3) implies the class of  
$\mathfrak{J}$ is
an extension of a class from $A(K)$ (recall the map $A(K)\rt A(K_{m-1})$ is 
just an injection $\z/p^a\z \hookrightarrow \z/p^{a+m-1}\z$). We let $\mathfrak{A}$ be a 
representative ideal of the class that extends to the class of $\mathfrak{J}$.

Since the
$p$-Hilbert class field of $K$ is contained in $\Kt$, and
the class of $\mathfrak{A}$, and therefore $\mathfrak{J}$, 
becomes principal in $\Kt$. The
extension $M_K/\Kt$ is also generated
by a $p^m$-th root of $x$, and the ideal $(x)$ in $\Kt$  
is now the $p^m$-th power of a
\textit{principal} ideal, 
\[(x)=(y)^{p^m}.\] The elements $x$ and $y^{p^m}$ then differ by
a unit , i.e. $x=uy^{p^m}$. But clearly the extension $M_K$ is 
also generated by the
$p^m$-th root of $x/y^{p^m}=u$, and so the field
$M_K$ is contained in $N_\infty$. $\Box$

\end{document}